\documentclass[ a4paper,11pt]{amsart}
\usepackage{amsmath}
\usepackage{amssymb}
\newcommand{\dis}{\displaystyle}

\begin{document}
\title{On Bernoulli Numbers and Its Properties}\texttt{}
\author{Cong Lin}
\address{\noindent{\small
Hwa Chong Junior College\\
661 Bukit Timah Road\\
Singapore 269734}}
\date{}
\maketitle

\begin{abstract}
In this survey paper, I first review the history of Bernoulli
numbers, then examine the modern definition of Bernoulli numbers
and the appearance of Bernoulli numbers in expansion of functions.
I revisit some properties of Bernoulli numbers and the history of
the computation of big Bernoulli numbers.
\end{abstract}

\section{Introduction}

Two thousand years ago, Greek mathematician Pythagoras first noted
about triangle numbers which is $1+2+3+\cdots+n$. Archimedes found
out \begin{equation}
1^2+2^2+3^2+\cdots+n^2=\frac{1}{6}n(n+1)(2n+1). \end{equation}
Later in the fifth century, Indian mathematician Aryabhata
proposed \begin{equation}
1^3+2^3+3^3+\cdots+n^3=[\frac{1}{2}n(n+1)]^2
\end{equation} which Jacobi gave the
first vigorous proof in 1834. It is not until five hundred years
later that Arabian mathematician Al-Khwarizm showed
\begin{equation}
1^4+2^4+3^4+\cdots+n^4=\frac{1}{30}n(n+1)(2n+1)(3n^2+3n-1).
\end{equation} Studies the more generalized formula for
\scriptsize $\dis \sum_{k=1}^{n-1}k^r$ \normalsize for any natural
number $r$ was only carried out in the last few centuries. Among
them,the investigation of Bernoulli numbers is much
significant.\\\\
In this paper, I present an elementary examination of the
development of Bernoulli numbers and a concise review of its
appearance in expansions of various functions . I also aim to
explore its properties and its applications in other fields of
mathematics \cite{Katz}.

\section{Bernoulli Numbers}

Swiss mathematician Jakob Bernoulli (1654-1705)once claimed that
instead of laboring for hours to get a sum of powers, he only used
several minutes to calculate sum of powers such as
$1^{10}+2^{10}+3^{10}+\cdots+1000^{10}
 = 91,409,924,241,424,243,424,241,924,242,500$. Obviously, he had used a
summation formula, knowing the first 10 Bernoulli numbers
\cite{Arts}.
\\\\It was already known to Jakob Bernoulli that
\begin{eqnarray}\sum_{k=1}^{n-1}k  =  \frac{1}{2}n(n-1)& = & \frac{1}{2} n^2 -
\frac{1}{2}n \\ \sum_{k=1}^{n-1}k^2 = \frac{1}{6}n(n-1)(2n-1)& = &
\frac{1}{3}n^3-\frac{1}{2}n^2+\frac{1}{6}n \\ \sum_{k=1}^{n-1}k^4
= \frac{1}{4}n^2(n-1)^2 & = & \frac{1}{4}n^4 - \frac{1}{2}n^3+
\frac{1}{4}n^2 \\
\sum_{k=1}^{n-1}k^4=\frac{1}{30}n(n-1)(2n-1)(3n^2-3n-1) & =
&\frac{1}{5}n^5-\frac{1}{2}n^4+{}\nonumber\\& & {}+\frac{1}{3}n^3
-\frac{1}{30}n \\ \sum_{k=1}^{n-1}k^5
=\frac{1}{12}n^2(2n^2-2n-1)(n-1)^2 & = & \frac{1}{6}n^6-
\frac{1}{2}n^5+{}\nonumber\\& & {}+\frac{5}{12}n^4-\frac{1}{12}n^2
\end{eqnarray} More generally, the $\dis \sum_{k=1}^{n-1}k^r $ can
be written in the form of
\begin{eqnarray} \sum_{k=1}^{n-1}k^r & = & \sum_{k=0}^{r} \frac{B_k}{k!}
\frac{r!}{(r-k+1)!} n^{r-k+1}{}\nonumber\\& = &
{}\frac{B_0}{0!}\frac{n^{r+1}}{r+1}+\frac{B_1}{1!}{n^r}+\frac{B_2}{2!}{r
n^{r-1}}+\dots +\frac{B_r}{r!}n \end{eqnarray} where the $B_k$
numbers which are independent of $r$ and called Bernoulli's
numbers. The first few Bernoulli numbers $B_n$ are
\begin{eqnarray} B_0& = & 1 \nonumber\\ B_1& = & -\frac{1}{2} \nonumber\\ B_2& =& \frac{1}{6} \nonumber\\ B_4& =&
-\frac{1}{30}  \nonumber\\  B_6& =& \frac{1}{42} \nonumber\\  B_8&
=& -\frac{1}{30}\nonumber
\end{eqnarray} \begin{eqnarray}
B_{10}& =& \frac{5}{66} \nonumber\\  B_{12}& =& -\frac{691}{2,730} \nonumber\\
B_{14}& =& \frac{7}{6}  \nonumber\\ B_{16}& =& -\frac{3,617}{510}  \nonumber\\
B_{18}& =& \frac{43,867}{798} \nonumber\\  B_{20}& =&
-\frac{17\nonumber4,611}{330} \nonumber\\ B_{22}& =&
\frac{854,513}{138}  \nonumber\\ B_{24}& =&
-\frac{236,364,091}{2,730}\nonumber\\ B_{26}& =& \frac{8,553,103}{6} \nonumber\\
B_{28}& =& -\frac{23,749,461,029}{870}  \nonumber\\ B_{30}& =&
\frac{8,615,841,276,005}{14,322} \nonumber\\B_{32}&
 = & -\frac{7,709,321,041,217}{510} \nonumber\\
B_{34} & = & \frac{2,577,687,858,367}{6} \nonumber
\\ B_{36}& = & -\frac{26,315,271,553,053,477,373}{1,919,190} \nonumber
\\B_{38} & = & \frac{2,929,993,913,841,559}{6} \nonumber
\\B_{40} & = & -\frac{261,082,718,496,449,122,051}{13,530} \nonumber \end{eqnarray} with
\begin{eqnarray} B_{2n+1}& = & 0 \end{eqnarray} for all positive
integer $n$. More numbers are given in \cite{Stegun}.

\section {\textbf{Expansion of Usual Functions}}

An equivalent definition of the Bernoulli's numbers is obtained
from the series expansion of the identity \begin{eqnarray}
\frac{x}{e^x-1} & \equiv & \dis \sum_{n=0}^\infty \frac{B_n
x^n}{n!}. \end{eqnarray} This leads to
\begin{eqnarray} \frac{x}{e^x-1} + \frac{x}{2} & = &
\frac{x}{2}(\frac{2}{e^x-1}+1) =\frac{x}{2}\frac{e^x+1}{e^x-1} =
\frac{x}{2}coth(\frac{x}{2}).
\end{eqnarray} It can be rewritten as \begin{equation}
\frac{x}{2}coth(\frac{x}{2}) = \sum_{n=0}^\infty \frac{B_{2n}
x^{2n}}{(2n)!}. \end{equation} If substitute $x$ with $2ix$, then
it gives \begin{equation} x cot (x) = \sum_{n=0}^\infty (-1)^n
\frac{B_{2n} (2x)^{2n}}{(2n)!}
\end{equation} for $x \in [-\pi,\pi]$. Thus the following expansions are
obtained, \begin{eqnarray} coth(x) = \sum_{n=0}^\infty
\frac{2B_{2n} (2x)^{2n-1}}{(2n)!}
\\ cot (x) =
\sum_{n=0}^\infty (-1)^n \frac{2 B_{2n} (2x)^{2n-1}}{(2n)!}
\end{eqnarray}
Now it's also possible to find the expansion for $tanh(x)$ and
$tan(x)$. As I observe that \begin{eqnarray} 2coth(2x)-coth(x) & =
& 2 \frac{cosh(2x)}{sinh(2x)}-\frac{cosh(z)}{sinh(z)} {}\nonumber
\\ & = & {} \frac{cosh^2(x)+sinh^2(x)}{sinh(x)cosh(x)}
-\frac{cosh(x)}
 {sinh(x)}
 =tanh(x)\nonumber , \end{eqnarray}
 I can have \begin{equation} tanh(x) =  \sum_{n=1}^\infty
\frac{2(4^n-1)B_{2n} (2x)^{2n-1}}{(2n)!},\end{equation} and then,
\begin{equation} tan(x)= \sum_{n=1}^\infty (-1)^n
\frac{2(1-4^n)B_{2n} (2x)^{2n-1}}{(2n)!} \end{equation} both for
$x \in (-\frac{\pi}{2}, \frac{\pi}{2})$ Bernoulli's numbers also
occur in the expansions of other classical functions such as
$csc(x)$,$csch(x)$,$ln|sin(x)|$,$ln|cos(x)|$,$ln|tan(x)|$,
$\frac{x}{sinh(x)}$, and etc.
\\ \\Another intriguing fact is that
zeta function $\zeta{(2k)}$ for any natural number $k$ can also be
expressed in Bernoulli numbers \cite{Borwein}
\begin{equation}\ \zeta{(2k)}=\sum_{n=1}^{\infty}\frac{1}{n^{2k}}=
\frac{4^k |B_{2k}| \pi^{2k}}{2(2k)!}\end{equation} However, there
is no similar expression known for $\zeta{(2k+1)}$. Rearrange the
equation, I get
\begin{equation} B_{2k}= \frac{(-1)^{k-1}2(2k)!}{(2 \pi)^{2k}}
\sum_{n=1}^{\infty} \frac{1}{n^{2k}}= \frac{(-1)^{k-1} 2 (2k)!}{(2
\pi^{2k})}\zeta{(2k)}.
\end{equation} This expression gives rise to an approximation of
Bernoulli numbers. When $k$ becomes large,
\begin{equation} \zeta{(2k)} \simeq 1 ,\end{equation} while with Stirling's formula
\begin{equation}(2k)!\simeq (2k)^{2k}e^{-2k}\sqrt{4\pi k}, \end{equation}
so I have \begin{equation} B_{2k}\simeq (-1)^{k-1}4(\frac{k}{\pi
e})^{2k}\sqrt{\pi k} .\end{equation}

\section{Properties of Bernoulli Numbers}

The Bernoulli numbers are given by the double sum
\begin{equation} B_n =\sum_{k=0}^n
\frac{1}{k+1}\sum_{r=0}^k(-1)^r\binom{k}{r} r^n
\end{equation}
where $\binom{k}{r}$ is a binomial coefficient. They also satisfy
the following interesting summations \cite{Lehmer, Carlitz}
\begin{eqnarray}
\sum_{k=0}^{n-1}\binom{n}{k} B_k & = & 0 ,
\\ \sum_{k=0}^{n}\binom{6n+3}{6k} B_{6k} & = & 2n+1 ,
\\ \sum_{k=0}^{n}\binom{6n+5}{6k+2} B_{6n+2} & = &
\frac{1}{3}(6n+5).
\end{eqnarray}
The famous Clausen-von Staudt's theorem regarding Bernoulli
numbers' fractional part was published by Karl von Staudt
(1798-1867) and Thomas Clausen (1801-1885) independently in 1840.
It allows to compute easily the fractional part of Bernoulli's
numbers and thus also permits to compute the denominator of those
numbers. It says, the value $B_{2k}$, added to the sum of the
inverse of prime numbers $p$ such that $(p-1)$ divides $2k$, is an
integer \cite{Hardy}. In other words,
\begin{equation} -B_{2k}\equiv \sum_{(p-1)|2k} \frac{1}{p} (mod 1)
\end{equation} For example, $k=8$,
\begin{eqnarray}-B_{16}& \equiv & \sum_{(p-1)|16} \frac{1}{p}\equiv \frac{1}{2}+\frac{1}{3}+
\frac{1}{5}+\frac{1}{17}\equiv \frac{47}{510} (mod 1){}\nonumber\\
&\Rightarrow & B_{16}\equiv \frac{463}{510} (mod 1). \nonumber
\end{eqnarray} One of the easy consequences of the Staudt's
theorem is that for every prime numbers $k$ of the form $3n+1$
\begin{equation}B_{2k} \equiv \frac{1}{6}(mod 1). \end{equation} This is because that $p-1$ divides
$2k=2(3n+1)$ only if $p-1$ is one of 1, 2, $3n+1$, $6n+2$, that is
p is one of 2, 3, $3n+2$, $6n+3$. But $6n+3$ is divisible by 3 and
$3n+2$ is divisible by 2 because 3n+1 is prime so the only primes
p candidates are 2 and 3. The first primes of the form 3n+1 are 7,
13, 19, 31, 37, 43, 61, 67, 73, 79, 97, 103... hence I have
\begin{eqnarray}
B_{14}& \equiv & B_{26}\equiv  B_{38}\equiv  B_{62}\equiv
B_{74}\equiv B_{86}\equiv B_{122} {}\nonumber\\ &\equiv &
B_{134}\equiv B_{146}\equiv B_{158}\equiv B_{194}\equiv
B_{206}\equiv \cdots \equiv \frac{1}{6} (mod 1)\nonumber
\end{eqnarray}
Staudt's theorem is very useful and significant in the sense that
it permits to compute exactly a Bernoulli's number as soon as
there is a sufficiently good approximation of it.
\\\\Besides its applications in series expansions, Bernoulli numbers
are also widely used in differential topology, mathematical
analysis, and number theory. Interestingly, it is also related to
the famous Fermat's last theorem.
\\\\Fermat's Last Theorem states
\begin{equation}
x^n + y^n  =  z^n \nonumber \end{equation} has no non-zero integer
solutions for $n>2$. Ever since Fermat expressed this result
around 1630, generations of mathematicians have dived
enthusiastically into the pursuit of a vigorous proof .
\\\\A breakthrough was made in 1850 by Ernst Kummer (1810-1893) when he
proved Fermat's theorem for \emph{n=}\emph{p}, whenever \emph{p}
is a regular prime. Kummer gave the beautiful regularity
criterion: \\\\\emph{p } is a regular prime if and only if\emph{
p} does not divide the numerator of $B_2,B_4,...,B_{p-3}$.
\\\\He showed that all primes before 37 where regular, hence Fermat's theorem was proved for those
primes. 37 is the first non regular prime because it divides the
numerator of
\begin{equation}
B_{32}=\frac{7709321041217}{510}=\frac{37¡Á208360028141
}{510}\nonumber
\end{equation}
The next irregular primes (less than 300) are
\begin{equation}59,67,101,103,131,149,157,233,257,263,271,283,293,...\nonumber
\end{equation}For example, 157 divides the numerators of $B_{62}$
and $B_{110}$.
\\\\By applying arithmetical properties of Bernoulli's numbers, Johann Ludwig
Jensen  (1859-1925) proved in 1915 that the number of irregular
primes is infinite \cite{last}.

\section{Computation of Bernoulli Numbers}
Bernoulli himself calculated the numbers up to $B_{10}$. Later,
Euler worked up to $B_{30}$. One century later, Adams made the
computa of all Bernoulli's numbers up to $B_{124}$ \cite{Adams}.
In 1996, Simon Plouffe and Greg J. Fee computed $B_{200,000}$, and
this huge number has about 800,000 digits. In July 10th 2002, they
improved the record to $B_{750,000}$ which has 3,391,993 digits by
a 21 hours computation on their personal computer \cite{Plouffe}.
The method is based on the relation between zeta function and
Bernoulli numbers, which allow a direct computation of the target
number without the need of calculating the previous numbers.


\begin{thebibliography}{99}
\bibitem{Katz} V.J. Katz, \emph{A History of Mathematics-An Introduction,} Addison-Wesley. (1998)
\bibitem{Arts} Jakob Bernoulli, \emph{Ars Conjectandi}, (1713)
\bibitem{Stegun} M. Abramowitz and I. Stegun, \emph{Handbook of Mathematical Functions,} Dover, New York, (1964)
\bibitem{Borwein} J.M. Borwein and P.B. Borwein, \emph{\emph{Pi and the AGM - A study in
Analytic Number Theory and Computational Complexity}}, A
Wiley-Interscience Publication, New York, (1987)
\bibitem{Lehmer} Lehmer, D. H. \emph{\emph{Lacunary Recurrences for the Bernoulli
Numbers, }}Ann. Math., vol. 36,p. 637-649 (1935)
\bibitem{Carlitz} Carlitz, L. \emph{Bernoulli Numbers,} Fib. Quart., vol. 6, p. 71-85 (1968)
\bibitem{Hardy} G.H. Hardy and E. M. Wright, An Introduction to the Theory of
Numbers, Oxford Science Publications, (1979)
\bibitem{last} P. Ribenboim, \emph{The new Book of Prime Number Records}, Springer, (1996)
\bibitem{Adams}J.C. Adams, \emph{On the calculation of Bernoulli's numbers up to B62 by
means of Staudt's theorem}, Rep. Brit. Ass., (1877)
\bibitem{Plouffe} Simon Plouffe, \emph{The 750,000'th Bernoulli Number as computed on July 10 2002,} http://www.lacim.uqam.ca/
\end{thebibliography}
\end{document}